\documentclass[10pt,a4paper,english]{article}
\usepackage[T1]{fontenc}
\usepackage[latin1]{inputenc}
\usepackage{amsthm}
\usepackage{amsmath}
\usepackage{amssymb}
\usepackage{amsfonts}
\usepackage[all]{xy}
\usepackage{babel}


\theoremstyle{plain}
\newtheorem{theorem}{Theorem}[section]
\newtheorem{axiom}[theorem]{Axiom}
\newtheorem{conjecture}[theorem]{Conjecture}
\newtheorem{corollary}[theorem]{Corollary}
\newtheorem{criterion}[theorem]{Criterion}
\newtheorem{lemma}[theorem]{Lemma}
\newtheorem{problem}[theorem]{Problem}
\newtheorem{proposition}[theorem]{Proposition}

\theoremstyle{definition}
\newtheorem{acknowledgement}[theorem]{Acknowledgement}
\newtheorem{algorithm}[theorem]{Algorithm}
\newtheorem{case}[theorem]{Case}
\newtheorem{claim}[theorem]{Claim}
\newtheorem{condition}[theorem]{Condition}
\newtheorem{conclusion}[theorem]{Conclusion}
\newtheorem{definition}[theorem]{Definition}
\newtheorem{example}[theorem]{Example}
\newtheorem{exercise}[theorem]{Exercise}
\newtheorem{notation}[theorem]{Notation}
\newtheorem{remark}[theorem]{Remark}
\newtheorem{solution}[theorem]{Solution}
\newtheorem{summary}[theorem]{Summary}

\newenvironment{prooof}[1][Proof]{\textbf{#1.} }
{\ \rule{0.5em}{0.5em}} 
\newenvironment{variant}[1][Variant]{\textbf{#1.} }
{\ \rule{0.5em}{0.5em}}

\DeclareMathOperator{\Alb}{Alb}

\DeclareMathOperator{\CH}{CH}

\DeclareMathOperator{\Cp}{Cp}

\DeclareMathOperator{\Decf}{\underline{Dec}}

\DeclareMathOperator{\Divf}{\underline{Div}}

\DeclareMathOperator{\FDivf}{\underline{FDiv}}

\DeclareMathOperator{\Extabk}{Ext_{\Abk}}

\DeclareMathOperator{\Hom}{Hom}

\DeclareMathOperator{\Homk}{Hom_{\fld}}

\DeclareMathOperator{\Homfabk}{\underline{Hom}_{\Abk}}

\DeclareMathOperator{\Lie}{Lie}

\DeclareMathOperator{\Pic}{Pic}

\DeclareMathOperator{\Picf}{\underline{Pic}}

\DeclareMathOperator{\Supp}{Supp}

\DeclareMathOperator{\Trace}{Tr}

\DeclareMathOperator{\chr}{char}

\DeclareMathOperator{\coker}{coker}

\DeclareMathOperator{\dv}{div}
\DeclareMathOperator{\fml}{fml}

\DeclareMathOperator{\im}{im}

\DeclareMathOperator{\rk}{rk}

\DeclareMathOperator{\val}{v}

\newcommand{\Alba}[1]{\Alb(#1)}
\newcommand{\Fm}[2]{\fmlG_{#1,#2}}

\newcommand{\AlbF}[1]{\Alb_{\fmlG}(#1)}

\newcommand{\Albsg}[1]{\Alb^{\mathrm{ESV}}(#1)}

\newcommand{\Divfc}[1]{\wh{\Divf_{#1}}}

\renewcommand{\H}{\mathrm{H}}

\newcommand{\Urq}[1]{\Alb(#1,#1_{\sing})}

\newcommand{\Urqlrpn}[1]{\Alb\lrpn{#1,(#1)_{\sing}}}

\newcommand{\Sppp}[1]{\mathrm{S} (#1)}

\newcommand{\CHSO}[1]{\CH_0 (#1,#1_{\sing})}

\newcommand{\Abk}{\mathcal{A} \mathit{b} / \fld}

\newcommand{\Mr}{\mathsf{Mr}}

\newcommand{\MrCH}[1]{\Mr^{\CH}}

\newcommand{\Nat}{\mathbb{N}}
\newcommand{\Zint}{\mathbb{Z}}

\newcommand{\Prj}{\mathbb{P}}

\newcommand{\dime}{\mathrm{d}}

\newcommand{\sing}{\mathrm{sing}}

\newcommand{\bdr}{\Phi}

\newcommand{\rmp}[1]{\phe^{#1}}

\newcommand{\fld}{\mathit{k}}

\newcommand{\Rfin}{R}

\newcommand{\Gp}{\mathrm{G}}

\newcommand{\fmlE}{\mathcal{E}}
\newcommand{\fmlEd}{\fmlE^{\vee}}

\newcommand{\fmlG}{\mathcal{F}}

\newcommand{\Ga}{\mathbb{G}_{\mathrm{a}}}
\newcommand{\Gac}{\widehat{\mathbb{G}}_{\mathrm{a}}}

\newcommand{\Gm}{\mathbb{G}_{\mathrm{m}}}
\newcommand{\GmS}[1]{\mathbb{G}_{\mathrm{m},#1}}

\newcommand{\Lin}{\mathbb{L}}

\newcommand{\Trs}{T}

\newcommand{\Upt}{U}

\newcommand{\Vcl}{V}

\newcommand{\Fsp}{f}
\newcommand{\fsp}{f_0}

\newcommand{\pnt}{q}
\newcommand{\pntt}{p}

\newcommand{\Es}{E}

\newcommand{\base}{\beta}
\newcommand{\Base}{B}

\newcommand{\Csp}{C_0}
\newcommand{\Cspt}{\widetilde{\Csp}}
\newcommand{\Ci}{C_i}
\newcommand{\Cit}{\widetilde{\Ci}}

\newcommand{\Crvv}{C}

\newcommand{\Crel}{\mathcal{Z}}
\newcommand{\Crll}{\mathcal{C}}
\newcommand{\CrX}{\Crll}
\newcommand{\CrY}{\Crel}

\newcommand{\linSys}{\left|\sL\right|}

\newcommand{\pt}{\widetilde{p}}

\newcommand{\Oc}{\widehat{\sO}}
\newcommand{\Ot}{\widetilde{\sO}}

\newcommand{\Xt}{\widetilde{X}}
\newcommand{\Ct}{\widetilde{C}}

\newcommand{\Cs}{C'}

\newcommand{\Galp}{\Gamma_{\alp}}
\newcommand{\Galpt}{\widetilde{\Gamma}_{\alp}}
\newcommand{\Gbet}{\Gamma_{\beta}}
\newcommand{\Gbett}{\widetilde{\Gamma}_{\beta}}

\newcommand{\alp}{\alpha}
\newcommand{\bet}{\beta}

\newcommand{\del}{\delta}

\newcommand{\tha}{\vartheta}

\newcommand{\lam}{\lambda}
\newcommand{\sig}{\sigma}

\newcommand{\phe}{\varphi}

\newcommand{\fm}{\mathfrak{m}}

\newcommand{\sD}{\mathcal{D}}

\newcommand{\sK}{\mathcal{K}}
\newcommand{\sL}{\mathcal{L}}

\newcommand{\sO}{\mathcal{O}}

\newcommand{\sQ}{\mathcal{Q}}

\newcommand{\bt}[1]{[#1]}

\newcommand{\bigbt}[1]{\big[#1\big]}

\newcommand{\pn}[1]{(#1)}
\newcommand{\lrpn}[1]{\left(#1\right)}
\newcommand{\bigpn}[1]{\big(#1\big)}
\newcommand{\Bigpn}[1]{\Big(#1\Big)}

\newcommand{\Biggpn}[1]{\Bigg(#1\Bigg)}

\newcommand{\st}[1]{\{#1\}}
\newcommand{\lrst}[1]{\left\{#1\right\}}
\newcommand{\bigst}[1]{\big\{#1\big\}}
\newcommand{\Bigst}[1]{\Big\{#1\Big\}}

\newcommand{\ra}{\rightarrow}

\newcommand{\dra}{\dashrightarrow}

\newcommand{\lra}{\longrightarrow}

\newcommand{\lmt}{\longmapsto}

\newcommand{\xra}[1]{\xrightarrow{#1}}

\newcommand{\wt}[1]{\widetilde{#1}}
\newcommand{\wh}[1]{\widehat{#1}}

\newcommand{\cut}{\cdot}
\newcommand{\isec}{\cap}
\newcommand{\dsum}{\bigoplus}

\newcommand{\tens}{\otimes}

\newcommand{\tms}{\times}

\newcommand{\lul}{?}                    
\newcommand{\llul}{\hspace{+.05em} ?}   
\newcommand{\lull}{? \hspace{+.05em}}   

\newcommand{\see}{\textrm{see }}
\newcommand{\seecite}{\textrm{see }}
\newcommand{\laurin}{.}                     
\newcommand{\laurink}{,}                   
\newcommand{\bDpl}{\;$}
\newcommand{\eDpl}{$\;}

\newcommand{\vs}{12pt}

\newcommand{\bThm}{\begin{theorem}}
\newcommand{\eThm}{\end{theorem}}
\newcommand{\bAck}{\begin{acknowledgement}}
\newcommand{\eAck}{\end{acknowledgement}}
\newcommand{\bAlg}{\begin{algorithm}}
\newcommand{\eAlg}{\end{algorithm}}
\newcommand{\bAxm}{\begin{axiom}}
\newcommand{\eAxm}{\end{axiom}}
\newcommand{\bCas}{\begin{case}}
\newcommand{\eCas}{\end{case}}
\newcommand{\bClm}{\begin{claim}}
\newcommand{\eClm}{\end{claim}}
\newcommand{\bCcl}{\begin{conclusion}}
\newcommand{\eCcl}{\end{conclusion}}
\newcommand{\bCdn}{\begin{condition}}
\newcommand{\eCdn}{\end{condition}}
\newcommand{\bCjc}{\begin{conjecture}}
\newcommand{\eCjc}{\end{conjecture}}
\newcommand{\bCor}{\begin{corollary}}
\newcommand{\eCor}{\end{corollary}}
\newcommand{\bCrt}{\begin{criterion}}
\newcommand{\eCrt}{\end{criterion}}
\newcommand{\bDef}{\begin{definition}}
\newcommand{\eDef}{\end{definition}}
\newcommand{\bExm}{\begin{example}}
\newcommand{\eExm}{\end{example}}
\newcommand{\bExc}{\begin{exercise}}
\newcommand{\eExc}{\end{exercise}}
\newcommand{\bLem}{\begin{lemma}}
\newcommand{\eLem}{\end{lemma}}
\newcommand{\bNot}{\begin{notation}}
\newcommand{\eNot}{\end{notation}}
\newcommand{\bPrb}{\begin{problem}}
\newcommand{\ePrb}{\end{problem}}
\newcommand{\bPrp}{\begin{proposition}}
\newcommand{\ePrp}{\end{proposition}}
\newcommand{\bRmk}{\begin{remark}}
\newcommand{\eRmk}{\end{remark}}
\newcommand{\bSol}{\begin{solution}}
\newcommand{\eSol}{\end{solution}}
\newcommand{\bSmr}{\begin{summary}}
\newcommand{\eSmr}{\end{summary}}
\newcommand{\bVar}{\begin{variant}}
\newcommand{\eVar}{\end{variant}}
\newcommand{\bPf }{\begin{prooof}}
\newcommand{\ePf }{\end{prooof}}

\begin{document}

\centerline{ } 
\vspace{35pt} 
\centerline{\LARGE{Description of Generalized}} 
\vspace{4pt} 
\centerline{\LARGE{Albanese Varieties by Curves}} 
\vspace{25pt}
\centerline{\large{Henrik Russell}} 
\vspace{30pt}


\begin{abstract}
Let $X$ be a projective variety over an algebraically closed base field, 
possibly singular. 
The aim of this paper is to show that 
the generalized Albanese variety $\Urq{X}$ of Esnault-Srinivas-Viehweg 
can be computed from one general curve $C$ in $X$, 
if the base field is of characteristic 0. 
We illustrate this by an example, 
which we also use to unravel some mysterious properties 
of $\Urq{X}$. 
\end{abstract}


\tableofcontents{} 
\newpage 

\setcounter{section}{-1}

\section{Introduction} 

In \cite{Ru} the author considered generalized Albanese varieties $\AlbF{X}$ 
associated with 
categories of rational maps from a variety $X$ to algebraic groups. 
If such a generalized Albanese variety 
is generated by a curve $\Crvv$, 
then the dimension of $\AlbF{X}$ is bounded by the dimension 
of the generalized Albanese of the curve $\Crvv$, 
which is easy to compute. 
For example if $\AlbF{X}$ is 
the classical Albanese $\Alba{X}$ of a smooth proper variety $X$ 
or the Albanese of Esnault-Srinivas-Viehweg $\Urq{X}$ 
of a (singular) projective variety $X$ 
($\seecite$\cite{ESV}, cf.\ \cite{Ru}), 
then $\Alba{C}$ resp.\ $\Urq{C}$ is isomorphic to the Picard variety $\Pic^0 C$ 
of $C$. 
In this way the existence of the classical Albanese was shown in \cite{La} 
and the existence of the Albanese of Esnault-Srinivas-Viehweg in \cite{ESV}. 

The purpose of the present paper is to show 
that the functorial description of generalized Albanese varieties from \cite{Ru} 
is not only a purely theoretical one, 
but allows a concrete computation, 
using the interplay with curves. 
Here we restrict ourselves to the case that 
the base field $k$ is algebraically closed of characteristic 0. 
The Albanese of Esnault-Srinivas-Viehweg $\Urq{X}$ is 
an extension of the classical Albanese $\Alba{\Xt}$, 
where $\Xt \ra X$ is a resolution of singularities, 
by an affine group $L_X$, 
whose Cartier dual $\Divf_{\Xt/X}^0$ is described in \cite[Prop.s~3.23, 3.24]{Ru}. 
The main result of this work is a significant simplification 
of the presentation of $\Divf_{\Xt/X}^0$ (Theorem~\ref{genCurve}). 
Moreover, the functorial description allows to explain 
some pathological properties of the Albanese of Esnault-Srinivas-Viehweg. 
This is accomplished in an example (Section~\ref{sec:Exm_Surf}).

\subsection{Leitfaden} 

\hspace{4.5mm} \textbf{Section~\ref{sec:Case-of-Curves}.} 
We recall some facts about the Picard variety of curves 
that allow us to compute 
$\Pic^0 C$ of a singular curve $C$. 
Here we decompose $\Pic^0 C$ as an extension 
of the Picard variety $\Pic^0 \Ct$ of the normalization $\Ct$ of $C$ 
by a linear group $L$ that takes care of the singularities of $C$. 

\textbf{Section~\ref{Reducing_to_Subvarieties}.} 
We consider a formal group $\fmlG \subset \Divf_X$ 
of relative Cartier divisors on $X$ 
and a curve $C$ in $X$. 
We give a sufficient condition for the injectivity of the restriction map 
$\fmlG \lra \Divf_C$ from $\fmlG \subset \Divf_X$ 
into the group sheaf of relative Cartier divisors on $C$ 
(Lemma~\ref{.C inj}). 

\textbf{Section~\ref{Computation}.} 
The Albanese of Esnault-Srinivas-Viehweg $\Urq{X}$ is 
the dual 
of the 1-motivic functor $\bigbt{\Divf_{\Xt/X}^0 \lra \Pic^0_{\Xt}}$, 
where $\Xt \ra X$ is a resolution of singularities, 
see \cite[Theorem 0.1]{Ru}. 
Here $\Divf_{\Xt/X}^0$ is the 
``kernel of the push-forward of relative Divisors from $\Xt$ to $X$'', 
if $X$ is a curve 
($\seecite$\cite[Proposition 3.23]{Ru}). 
For higher dimensional $X$ the definition of $\Divf_{\Xt/X}^0$ 
is derived from the one for curves 
by intersecting the formal groups $\Divf_{\Ct/C}^0$ 
associated with curves $C$ in $X$, 
where the intersection ranges over all Cartier curves $C$ in $X$ 
relative to the singular locus of $X$ 
($\seecite$\cite[Proposition 3.24]{Ru}). 
So a priori this object looks hard to grasp. 
We explain how $\Divf_{\Xt/X}^0$ 
can be computed from one single complete intersection curve in $X$ 
(Corollary~\ref{oneCurve}). 
Moreover, the curves with this property are dense 
in any space of sufficiently ample complete intersection curves 
(Theorem~\ref{genCurve}). 

Sections \ref{sec:Case-of-Curves}, \ref{Reducing_to_Subvarieties} 
and \ref{Computation} provide the necessary tools in order to reduce 
the computation of $\Urq{X}$ to the curve case. 
This is demonstrated in an example in the next section. 

\textbf{Section~\ref{sec:Exm_Surf}.} 
The classical Albanese of smooth projective varieties $X_i$ 
is compatible with products, i.e.\ 
$\Alb\lrpn{\prod X_i} = \prod \Alba{X_i}$. 
More generally, all universal objects 
for categories of rational maps to \emph{semi-abelian varieties} 
have this property. 
However, due to additive subgroups it is possible 
for singular projective varieties $X_i$ that 
\;$\dim \Urqlrpn{\prod X_i} > \sum \dim \Urqlrpn{X_i}$. 
Moreover, if $X$ is a smooth projective variety of dimension $d$ 
and $\sL$ a very ample line bundle on $X$, 
then for a complete intersection $C$ 
of $d - 1$ general divisors in the linear system $|\sL|$ 
the Gysin map \;$\Alba{C} \lra \Alba{X}$\; will be surjective. 
This is not true in general for the Albanese of Esnault-Srinivas-Viehweg 
of a singular $X$, 
but a sufficiently high power $\sL^{\tens N}$ will again have this property. 

We discuss an example 
that illustrates these pathological properties. 
This example was computed in the diploma of Alexander Schwarzhaupt 
\cite{Sch} by means of the Hodge theoretic description given in \cite{ESV}. 
The functorial description from \cite{Ru} 
yields a considerable simplification for the computation 
and gives an explanation for the strange behaviour of $\Urq{X}$. 
In particular we obtain a formula for 
the dimension of $\Urqlrpn{\prod X_i}$ 
as a function of the dimensions of $\Urqlrpn{X_i}$ 
(Proposition \ref{dimAlb(CxC)}). 
Using a formula from \cite{ESV}, 
Schwarzhaupt computes in this example bounds $N_0$, $N_1$ 
such that the Gysin map \;$\Urq{C} \lra \Urq{X}$\; 
is surjective for any complete intersection curve $C$ 
of general divisors in $|\sL^{\tens N}|$ if $N > N_1$ 
and not surjective if $N < N_0$. 
With our method we can optimize the bound $N_1$, 
in fact we can show that the Gysin map is surjective if $N \geq N_0$, 
i.e.\ the condition we obtain is necessary and sufficient for the 
surjectivity of the Gysin map (Proposition \ref{Gysin surjective}).

\vspace{\vs} 

\textbf{Acknowledgement.} 
This paper ties in with my PhD thesis. 
I thank again H\'el\`ene Esnault and Eckart Viehweg 
wholeheartedly for their support and guidance. 
In particular I use this article to commemorate Eckart Viehweg, 
who suggested this point of view to me. 
His brilliant work bears fruit, 
and in this sense he is still with us and helping us.


\section{Picard Variety of Curves} 
\label{sec:Case-of-Curves}

Let $C$ be a projective curve over a 
field $k$. 
The fact that divisors equal 0-cycles on $C$ yields an identification 
of the Picard scheme $\Pic C$ with the relative Chow group $\CHSO{C}$ 
from \cite{LW}. 

We describe the algebraic group $\Pic^0C$ as an 
extension of the abelian variety $\Pic^0\Ct$ 
by an affine algebraic group $L$, 
using an intermediate curve $\Cs$ (the semi-normalization) 
between $C$ and its normalization $\Ct$. 
The methods of this Section 
are taken from \cite[Section~9.2]{BLR}.

\subsection{Normalization and Semi-Normalization of a Curve} 
\label{sub:Normalzation}

Let $C$ be a curve over a perfect field. 

\bDef 
\label{ord-mult-pt}
A point $p$ of a curve $C$ is called an \emph{ordinary multiple point}, 
if it marks a transversal crossing of smooth formal local branches. 
More precisely, $p\in C$ is an \emph{ordinary $m$-ple point} if 
\bDpl \widehat{\sO}_{C,p} \;\cong\; 
k\bt{\bt{t_{1},\ldots,t_{m}}} \big/ \sum_{i\neq j}\lrpn{t_{i} \,t_{j}} \laurin 
\eDpl 
\eDef 

\bDef 
\label{semi-normalization}
The normalization $\nu: \Ct \ra C$ 
factors as $\Ct \xra{\sig} \Cs \xra{\rho} C$ 
for a unique curve $\Cs$ 
which is homeomorphic to $C$ and has only ordinary multiple points 
as singularities, see \cite[Section 9.2, p.~247]{BLR}. 
$\Cs$ is called the \emph{largest curve homeomorphic to $C$}, 
or the \emph{semi-normalization} of $C$. 
\eDef 

Since $C$ is reduced, the smooth locus is dense in $C$, hence the
singular locus $S$ is finite. The curve $\Cs$ is obtained from $\Ct$
by identifying the points $\pt_{i}\in\Ct$ lying over $p\in S$. 
For an explicit description of $\Cs$ see \cite[Section 9.2, p.~247]{BLR}. 

\bNot 
We write 
\bDpl \sO := \sO_C  \laurink \hspace{1mm} 
   \sO' := \rho_*\sO_{\Cs}  \laurink \hspace{1mm} 
   \Ot := \nu_*\sO_{\Ct}  \laurin 
\eDpl 
\eNot

\subsection{Decomposition of the Picard Variety} 
\label{sub: Ext}

\bPrp 
\label{multi pt} 
Let $\Cs$ be a connected projective curve having
only ordinary multiple points as singularities. 
Let $\Ct$ be the normalization of $\Cs$. 
Then $\Pic^0 \Cs$ is an extension of the abelian variety $\Pic^0\Ct$ 
by a torus $\Trs$. 
If $k$ is algebraically closed, 
$\Trs \cong (\Gm)^t$ is a split torus of rank 
\bDpl t = \sum_{m\geq 1} (m-1)\,\# S_m - \#\Cp(\Cs) + 1 \laurink 
\eDpl 
where $\Cp(\Cs)$ is the set of irreducible components of $\Cs$ 
and $S_m$ is the set of $m$-ple points 
($\see$Definition~\ref{ord-mult-pt}). 
\ePrp 

\bPf 
(See also \cite[Section~9.2, Proposition~10]{BLR} for the first statement.) 

\noindent 
If $S$ denotes the singular locus of $\Cs$ 
and $\Cp\pn{\Ct}$ the set of components of $\Ct$, 
we obtain from 
\;$1 \ra \sO'^* \ra \Ot^* \ra \sQ^* \ra 1$\; 
with $\sQ^* = \prod_{p \in S} \bigpn{\Ot_{p}}^* \big/ \lrpn{\sO'_{p}}^*$ 
the long exact cohomology sequence 
\[ 1 \lra k_{\Cs}^* \lra \prod_{Z\in\Cp\pn{\Ct}}k_{Z}^* \lra 
   \prod_{p\in S}\Trs_{p}\pn{k} \lra \Pic\pn{\Cs} \lra 
   \prod_{Z\in\Cp\pn{\Ct}}\Pic\pn{Z} 
   \lra 1 
\] 
where $k_{\Cs}^* = \H^0\pn{\Cs,\sO_{\Cs}^*}$, 
$k_{Z}^* = \H^0\pn{Z,\sO_{Z}^*}$, and 
\[ \Trs_{p}\pn{k}  =  \bigpn{\Ot_{p}}^* \big/ \lrpn{\sO'_{p}}^* \\ 
  =  \frac{\prod_{q\ra p}k(q)^*}{k(p)^*} \\ 
  \cong  \lrpn{k^*}^{m_{p}-1}
\] 
since each $p\in S$ is an ordinary multiple point 
($\see$Definition~\ref{ord-mult-pt}). 
Here $\prod_{p\in S}\Trs_{p}\pn{k}$ maps to 
the connected component of the identity of $\Pic\pn{\Cs}$. 
Then the affine part $T$ of $\Pic^0\Cs$ 
is the torus given by 
\[ \Trs\pn{k} 
  =  \coker\Biggpn{\prod_{Z\in\Cp\pn{\Ct}}k_{Z}^* \lra
            \prod_{p\in S}\Trs_{p}\pn{k}} \\ 
  =  \frac{\prod_{p\in S}\Trs_{p}\pn{k}}
            {\prod_{Z\in\Cp\pn{\Ct}}k_{Z}^* \Big/ k_{\Cs}^*} \\ 
  \cong  \lrpn{k^*}^{t} 
\] 
with \;$t \,=\, \sum_{m\geq 1} (m-1)\,\# S_m - \#\Cp(\Ct) + 1$\; 
and \;$\#\Cp(\Ct) = \#\Cp(\Cs)$. 
\ePf 

\bPrp 
\label{cusp} 
Let $C$ be a projective curve,
let $\Cs$ be the largest homeomorphic curve 
between $C$ and its normalization $\Ct$. 
Then $\Pic^0 C$ is an extension of $\Pic^0 \Cs$ 
by a unipotent group $\Upt$. 
If $k$ is algebraically closed, 
$\Upt$ is characterized by 
\bDpl \Upt\pn{k} = \prod_{p\in S} \pn{1+\fm_{\Cs,p}} / \pn{1+\fm_{C,p}} \laurink 
\eDpl 
if moreover $\chr(k) = 0$, the exponential map yields an isomorphism 
\bDpl \Upt\pn{k} \cong \prod_{p \in S} {\fm_{\Cs,p}} / {\fm_{C,p}} \laurink 
\eDpl 
where $S = C_{\sing}$ is the singular locus of $C$. 
\ePrp 

\bPf 
(See also \cite[Section~9.2, Proposition~9]{BLR} for the first statement.) 

\noindent 
The exact sequence 
\;$ 1 \ra \sO^* \ra {\sO'}^* \ra \sQ^* \ra 1 $ 
with $\sQ^* = \prod_{p \in S} {\sO'_p}^* \big/ {\sO_p}^*$ 
yields the exact cohomology sequence 
\bDpl 1 \ra \H^0\pn{\sQ^*} \ra \H^1\pn{\sO^*} \ra 
   \H^1\pn{\sO'^*} \ra1 \laurink
\eDpl 
since $\rho:\Cs\lra C$ is a homeomorphism, thus 
$\H^0\bigpn{\sO'^*} \cong \H^0\pn{\sO^*}$. 
Moreover it holds 
\bDpl \H^0\pn{\sQ^*} = \prod_{p\in S} {\pn{\sO'_p}^*} / {\sO_p^*} 
   = \prod_{p\in S} \pn{1+\fm_{\Cs,p}} / \pn{1+\fm_{C,p}} \laurin 
\eDpl 
As $\Upt\pn{k} := \H^0\pn{\sQ^*}$ is a connected unipotent group, 
the image of $\Upt\pn{k}$ is contained in the connected component 
of the identity of $\Pic(C) = \H^1\pn{\sO^*}$. 
\ePf

\bPrp 
\label{Main} 
Let $C$ be a projective curve and $\Ct$ its normalization. 
Then the Picard variety $\Pic^0 C$ is an extension 
of the abelian variety $\Pic^0\Ct$ 
by an affine algebraic group 
$L = \Trs \tms_k \Upt$, 
which is the product of a torus $\Trs$ and a unipotent group $\Upt$. 
If $k$ is algebraically closed, 
$\Trs$ and $\Upt$ are characterized as in Propositions~\ref{multi pt} and \ref{cusp}. 
\ePrp 

\bPf 
Follows directly from Propositions~\ref{multi pt} and \ref{cusp}, 
cf.\ \cite[Section~9.2, Corollary~11]{BLR}. 
%
\ePf 

\bRmk 
\label{Urq(C)}
Let $C$ be a projective curve. 
Then there is a canonical isomorphism $\Pic^0 C \cong \Urq{C}$, 
$\seecite$\cite[Introduction]{ESV}. 
In particular, the affine part $L$ of $\Pic^0 C$ is Cartier dual to $\Divf_{\Ct/C}^0$, 
by \cite[Thm.~0.1]{Ru} for $\chr(k) = 0$ and by \cite[Thm.~0.1]{Ru3} in general. 
\eRmk

\section{Restriction to Curves} 
\label{Reducing_to_Subvarieties}

Let $X$ be a regular projective variety of dimension $d$ 
over an algebraically closed field $k$ of characteristic 0. 

Let $V$ be a subvariety of $X$. 
Let $\Divf_X$ be the functor of relative Cartier divisors 
as described in \cite[No.~2.1]{Ru}. 
There is a canonical class map $\Divf_X \lra \Picf_X$, 
set $\Divf_X^0 := \Divf_X \tms_{\Picf_X} \Picf_X^0$. 

$\Decf_{X,V}$ is the subfunctor of $\Divf_X$ 
consisting of those families of Cartier divisors 
whose support intersects $V$ properly, 
i.e.\ $\Supp\pn{\sD}$ ($\seecite$\cite[Def.~2.10]{Ru}) 
does not contain any associated point of $V$ 
for all $\sD\in\Decf_{X,V}\pn{\Rfin}$, $\Rfin$ a finite dimensional $k$-algebra 
($\seecite$\cite[Def.~2.11]{Ru}). 
$\lul\cdot V:\Decf_{X,V}\lra\Divf_{V}$ is the pull-back of relative 
Cartier divisors from $X$ to $V$ 
($\seecite$\cite[Def.~2.12]{Ru}). 

Let $\delta\in\Lie\lrpn{ \Divf_X} = \Gamma \lrpn{\sK_X/\sO_X}$ 
be a deformation of the zero divisor in $X$. 
Then $\delta$ determines an effective divisor by the poles of its local sections. 
Hence for each generic point $\eta$ of height 1 in $X$, 
with associated discrete valuation $\val_{\eta}$, 
the expression $\val_{\eta}(\delta)$ is well defined and 
$\val_{\eta}(\delta)\leq 0$. Thus we obtain a homomorphism 
$\val_{\eta}:\Lie\lrpn{ \Divf_X} \lra \Zint$. 

\bDef 
\label{|sL|_S} 
For an ample line bundle $\sL$ on $X$ and
an integer $c$ with $1\leq c\leq\dim X$ write 
\bDpl 
\left|\sL\right|^{c}=\Prj\lrpn{\H^0\lrpn{X,\sL}}\times\ldots
\times\Prj\lrpn{\H^0\lrpn{X,\sL}} 
\hspace{2mm} \textrm{($c$ copies)} \laurin 
\eDpl 
Let $H_{1},\ldots,H_{c}\in\left|\sL\right|$ and $V=\bigcap_{i=1}^{c}H_{i}$. 
By abuse of notation we write $V\in|\sL|^{c}$ instead
of $\lrpn{H_{1},\ldots,H_{c}} \in |\sL|^{c}$. 
\eDef  

%
%

\bLem 
\label{.C inj}
Let $\fmlG$ be a 
formal subgroup of $\Divf_X^0$ 
s.t.\ $\fmlG \cong \Zint^t \tms_k \pn{\Gac}^v$ for $t,v \in \Nat$, 
where $\Gac$ denotes the completion of $\Ga$ at $0$. 
Let $S$ be the set of generic points of $\Supp(\fmlG)$ and $S_{\inf}$ 
the corresponding set for $\Supp(\fmlG_{\inf})$. If $\eta$ is a generic 
point of height 1 in $X$, denote by $E_{\eta}$ the associated prime divisor. 
For an ample line bundle $\sL$ on $X$ 
there is an open dense $U \subset |\sL|^{d-1}$ such that 
\;$\lrpn{\lul\cut C}|_{\fmlG}:\fmlG\lra\Divf_C$\; is injective for $C \in U$ if 
\[ \# \lrpn{ C \cap E_{\eta} } \, \geq \, 
   \dim_k \, \pn{ \Lie \fmlG }_{\eta}^{-\nu} 
\] 
for all $\eta \in S_{\inf}$ and all $-\nu \in \val_{\eta} (\Lie \fmlG)$, 
where $\pn{\Lie \fmlG}_{\eta}$ is the image of the localization 
\;$\Lie \fmlG \subset \Gamma \lrpn{\sK_X/\sO_X} \lra \lrpn{\sK_X/\sO_X}_{\eta}$, 
$\delta \lmt [\delta]_{\eta}$\; 
at the height 1 point $\eta$, and 
$\pn{\Lie \fmlG}_{\eta}^{-\nu} = 
 \bigpn{ \pn{\Lie \fmlG}_{\eta} \cap \fm_{X,\eta}^{-\nu}\big/\sO_{X,\eta} }  \big/ 
 \bigpn{ \pn{\Lie \fmlG}_{\eta} \cap \fm_{X,\eta}^{-\nu+1}\big/\sO_{X,\eta} }$. 

\noindent 
Here 
\;$\lrpn{\sK_X/\sO_X}_{\eta} = \bigcup_{\nu>0} \fm_{X,\eta}^{-\nu}\big/\sO_{X,\eta}$\; 
and \;$\fm_{X,\eta}^{-\nu} = \st{f \in \sK_{X,\eta} \;|\; \val_{\eta}(f) \geq -\nu}$. 
\eLem 

\bPf  
Let $C = \bigcap_{i=1}^{d-1} H_i$ 
be a complete intersection curve in $|\sL|^{d-1}$. 
As an ample divisor $H_i$ intersects each closed subscheme of codimension
1 and $H_{i}$ restricted to $H_{1}\cap\ldots\cap H_{i-1}$ is again
ample for all $i=2,\ldots,d-1$, it follows by induction that
$C\cap\Supp\pn{\sD}\neq\varnothing$ for all $0 \neq \sD \in \fmlG(\Rfin)$, 
$\Rfin$ a finite dimensional $k$-algebra. 
If the intersection points of $C$ and $\Supp(D)$ are in general position 
for each $0 \neq D \in \fmlG(k)$, 
then $\lrpn{\lul\cut C}|_{\fmlG(k)}:\fmlG(k)\lra\Divf_C(k)$ 
is injective. 

For the infinitesimal part of $\fmlG$ consider the following diagram: 
\[ \xymatrix{ \Gamma \lrpn{\sK_X/\sO_X}
   \supset \Lie\fmlG \ar[r]^-{\lul\cut \, C} \ar@<+7.5ex>[d]^{\wr} & 
   \Gamma \lrpn{\sK_C/\sO_C} \ar[d]^{\wr}  \\  
   \dsum_{\eta \in S}\lrpn{\sK_X/\sO_X}_{\eta} 
   \supset \im \lrpn{ \Lie \fmlG } \ar[r] &  
   \dsum_{q \in C} \lrpn{\sK_C/\sO_C}_q \laurin } 
\]
For each $\eta \in S_{\inf}$ choose a local parameter $t_{\eta}$ of 
$\fm_{X,\eta}$. 
Since $C$ intersects $E_{\eta}$ properly, 
and if $C$ intersects $E_{\eta}$ transversally in general points, 
we may assume that each $q \in C \cap E_{\eta}$ is a regular 
closed point of $C$. Then 
the image $t_q \in \sO_C$ of $t_{\eta} 
\in \sO_X$ is a local parameter of $\fm_{C,q}$ and 
$\val_{q} (\Lie \fmlG\cut C) = \val_{\eta} (\Lie \fmlG)$. 
Set $-n_{\eta} = \min\lrst{\val_{\eta}(\Lie\fmlG)}$. 
We may consider $\Lie \fmlG \cut C$ as a $k$-linear subspace of the 
$k$-vector space 
\[  \dsum_{\eta \in S_{\inf}} \dsum_{q \in C \cap E_{\eta}} 
    \left.t_q^{-n_{\eta}}\sO_{C,q}\right/\sO_{C,q} \laurin 
\] 
Then the map 
\begin{eqnarray*}
\im \lrpn{ \Lie \fmlG \lra \lrpn{\sK_X/\sO_X}_{\eta} } 
  & \lra & 
  \dsum_{q \in C \cap E_{\eta}} t_q^{-n_{\eta}}\sO_{C,q}/\sO_{C,q} \\ 
f \, t_{\eta}^{-\nu} & \lmt & \sum_{q \in C \cap E_{\eta}} [f]_q \, t_q^{-\nu}
\end{eqnarray*} 
is injective if 
\bDpl \dim_k \, \pn{ \Lie \fmlG }_{\eta}^{-\nu} 
   \,\leq\,  \# \lrpn{ C \cap E_{\eta} } 
\eDpl 
for all $-\nu \in \val_{\eta}\pn{\Lie \fmlG}$ 
and the intersection points of $C$ and $E_{\eta}$ are in general position. 
%
$\lrpn{\lul\cut C}|_{\fmlG}: \Lie\fmlG \lra \Lie\lrpn{\Divf_C}$ 
is injective if these maps are injective for all $\eta \in S_{\inf}$. 
Since $\fmlG \cong \Zint^t \tms_k \pn{\Gac}^v$ by assumption, 
$\fmlG$ is already determined by $\fmlG(k)$ and $\Lie\fmlG$. 
\ePf

\section{Computation of the Kernel of the Push-forward of Divisors} 
\label{Computation}

Let $X$ be a projective variety of dimension $d$ 
over an algebraically closed field $k$ of characteristic 0. 
Let $S$ denote the singular locus of $X$. 
For a line bundle $\sL$ on $X$ we define 
\bDpl \linSys_S = \bigst{H\in\linSys \;\big|\; H \textrm{ intersects } 
                                      S \textrm{ properly} } \laurin 
\eDpl 
Let $\pi:\Xt\lra X$ be a projective resolution of singularities. 
For a curve $C$ 
we denote the normalization of $C$ by \;$\nu_C: \Ct \lra C$. 
The functor of \emph{formal divisors} on $C$ is the formal group given by 
\bDpl \FDivf_{C} = \dsum_{\pntt\in C(k)} 
   \Homfabk\bigpn{\Oc_{C,\pntt}^*,k^*} 
\eDpl 
($\seecite$\cite[Def.~2.1]{Ru3}). 
A finite morphism 
$\zeta: Z \lra C$ induces 
an obvious \emph{push-forward} of formal divisors 
$\zeta_*: \FDivf_Z \lra \FDivf_C$ 
($\seecite$\cite[Def.~2.4]{Ru3}). 
If $C$ is normal, 
there is a canonical homomorphism 
\;$\fml: \Divfc{C} \lra \FDivf_C$\; 
given by \;$\sD  \lmt  \sum_{\pntt\in C(k)} (\sD,\lull)_{\pntt}$, 
where $(\llul,\lull)_{\pntt}$ is the local symbol at $\pntt \in C$ 
($\seecite$\cite[Prop.~2.5]{Ru3}). 

We use the definition of $\Divf_{\Xt/X}^0$ given in 
\cite[Def.~2.6, 2.7]{Ru3}, 
which coincides with the one given in 
\cite[Prop.~3.23, 3.24]{Ru}, 
as is easily verified 
($\seecite$\cite[Rmk.~2.8]{Ru3}): 

\bDef 
\label{Div_Y/X}
If $C$ is a projective curve, then 
\[ \Divf_{\Ct/C}^{0} =  
   \ker\lrpn{\Divf_{\Ct}^{0} \overset{\fml} \lra
       \FDivf_{\Ct} \overset{\nu_{*}} \lra \FDivf_{C}} 
\] 
where \;$\nu: \Ct \ra C$\; is the normalization. 
For higher dimensional $X$ 
\[  \Divf_{\Xt/X}^0 =
    \bigcap_{C} \lrpn{ \llul\cut\Ct }^{-1} \Divf_{\Ct/C}^0 
\] 
where the intersection ranges over all Cartier curves in $X$ relative to 
the singular locus of $X$ ($\seecite$\cite[Definition 3.1]{Ru}), 
and $\bigpn{\llul \cut \Ct}$ is the pull-back of relative Cartier divisors on $\Xt$ 
to those on $\Ct$. 
\eDef 

The functor $\Divf_{\Xt/X}^0$ is a torsion-free dual-algebraic 
(\cite[Def.~1.20]{Ru2}) formal group 
($\seecite$\cite[Thm.~4.5]{Ru3} or \cite[Prop.~3.24]{Ru}), 
this means that the Cartier dual of $\Divf_{\Xt/X}^0$ 
is a connected algebraic affine group 
(by definition and \cite[(5.2)]{L}). 
Equivalently, $\Divf_{\Xt/X}^0 \cong \Zint^t \tms \pn{\Gac}^v$ 
for some $t,v \in \Nat$ (cf.\ \cite[(4.2)]{L}). 

As $\chr(k) = 0$, a formal group $\fmlE$ is completely determined 
by its $k$-valued points $\fmlE(k)$ and its Lie-algebra $\Lie\pn{\fmlE}$ 
($\seecite$\cite[Cor.~1.7]{Ru}). 
If $\fmlE$ is torsion-free and dual-algebraic, 
then $\fmlE(k)$ is a free abelian group of finite rank 
and $\Lie\pn{\fmlE}$ is a finite dimensional $k$-vector space. 
Then the dimension of the Cartier dual $\fmlEd$ of $\fmlE$ is 
\;$\dim \fmlEd = \rk \fmlE(k) + \dim \Lie\pn{\fmlE}$ 
(cf.\ \cite[(5.2)]{L}). 

\bPrp[Bertini's Theorem] 
\label{BertiniThm} 
Let $\sL$ be a line bundle on $X$. 
Then for almost all $C\in\linSys$ the inverse image 
\bDpl C_{\Xt} = C \tms_X \Xt = \pi^{-1}C \in \left|\pi^*\sL\right| 
\eDpl 
is smooth. 
\ePrp 

\bPf 
\cite[II, Theorem~8.18]{H} and induction. 
\ePf 

%

\bPrp 
\label{openProperty}
Let $\Base$ be a variety parametrizing Cartier curves in $X$, 
$\Csp\in\Base$ and $\sD\in\Divf^0_{\Xt}$ such that 
$\bigpn{\nu_{\Csp,*} \circ \fml} \bigpn{\sD\cut\Cspt} \neq 0$. 
Then there exists an open neighbourhood $U\ni\Csp$ in $\Base$ 
such that 
$\bigpn{\nu_{C,*} \circ \fml} \bigpn{\sD\cut\Ct} \neq 0$ for all $C\in U$. 
In other words: 
The zero locus of the function 
\bDpl \bigpn{\nu_{\lul,*}\circ \fml} \bigpn{\sD \cut \wt{\lul}} : 
   \Base \lra \dsum_{C\in \Base} \FDivf_C 
\eDpl 
is closed. 
\ePrp 

\bPf 
Let $\CrX \overset{\base} \lra \Base$ be the universal curve over $\Base$, 
let $\CrY := \CrX \tms_X \Xt$ 
be the relative curve over $\Base$ of preimages of $\CrX$ in $\Xt$. 
We may assume that the fibre $\CrY_{\base\pn{\Csp}}$ of $\CrY$ 
over the point $\base\pn{\Csp} \in \Base$ corresponding to $\Csp$ 
is normal, i.e. $\Cspt = \CrY_{\base\pn{\Csp}}$. 
Otherwise this can be achieved by blowing up $\Xt$. 

If $C$ is a curve, $S \subset C(k)$ and $F \in \FDivf_C$, we denote by 
$F_S := \sum_{\pntt \in S} F_{\pntt}$ 
$\in \Homfabk\bigpn{\Oc_{C,S}^*, k^*}$ 
the sum of those components 
$F_{\pntt} \in \Homfabk\bigpn{\Oc_{C,\pntt}^*, k^*}$ 
with $\pntt \in S$. 

By definition, 
$\bigpn{\nu_{\Csp,*} \circ \fml} \bigpn{\sD\cut\Cspt} \neq 0$ 
implies that there is $\pntt \in \Csp(k)$ and $\fsp \in \sO_{\Csp,\pntt}^*$ such that 
$\sum_{\pnt \ra \pntt} \bigpn{\sD\cut\Cspt,\fsp}_{\pnt} \neq 0$, 
where $\pnt \ra \pntt$ are the points $\pnt \in \Cspt$ over $\pntt \in \Csp$. 
Let $\Sppp{\sD}$ denote the support of $\sD$. 
Then necessarily $\pnt \in \Sppp{\sD}$. 
Let $\Lin_R := \Gm \lrpn{\llul \tens R}$ be the Weil restriction of $\GmS{R}$ 
from $R$ to $k$. 
We are going to construct a regular map $\bdr: T \lra \Lin_R$ 
from a neighbourhood $T$ of $\base\pn{\Csp}$ in $\Base$ to $\Lin_R$, 
such that 
\bDpl \bdr(b) = 
   \bigpn{\pn{\nu_{\CrX_b,*} \circ \fml} 
   \pn{\sD\cut\CrY_b}}_{\pi \Sppp{\sD} \,\cap\, \CrX_b} 
   \pn{\Fsp|_{\CrX_b}} 
\eDpl 
for all $b \in T$ for some $\Fsp \in \sO_{\CrX_T,\pi \Sppp{\sD}_{\CrX}}^*$, 
and $\bdr\bigpn{\base\pn{\Csp}} = 
\bigpn{\pn{\nu_{\Csp,*} \circ \fml} \pn{\sD\cut\Cspt}}_{\pntt} \pn{\fsp}$. 
Then $\bdr\bigpn{\base\pn{\Csp}} \neq 0$ 
implies that there is an open neighbourhood $U \ni \base\pn{\Csp}$ in $\Base$ 
such that $\bdr(u) \neq 0$ for all $u \in U$, 
proving the assertion. 

Let $\Gp\pn{\sD} \in \Extabk\pn{\Alba{\Xt}, \Lin_R} \cong \Pic_{\Alba{\Xt}}^0(R) 
\cong \Pic_{\Xt}^0(R)$ 
be the algebraic group corresponding to $\sO_{\Xt \tens R}(\sD)$. 
The canonical 1-section of $\sO_{\Xt \tens R}(\sD)$ induces a rational map 
$\rmp{\sD}: \Xt \dra \Gp\pn{\sD}$, 
which is regular away from $\Sppp{\sD}$. 
If $f \in \sO_{\CrY_b,\pnt}^*$ for some $b \in \Base$ and $\pnt \in \CrY_b$, 
then according to \cite[Lem.~3.16]{Ru} 
$\pn{\rmp{\sD}|_{\CrY_b}, f}_{\pnt}$ is contained 
in the fibre of $\Gp\pn{\sD}$ over $0 \in \Alba{\Xt}$, which is $\Lin_R$, 
and $\pn{\rmp{\sD}|_{\CrY_b}, f}_{\pnt} = \pn{\sD \cut \CrY_b, f}_{\pnt}$. 

Let $M$ be a modulus for the rational map $\rmp{\sD}|_{\Cspt}$. 
By the Approximation Lemma we find $f_p \in \sO_{\Csp,\pntt}^*$ 
with $\fsp / f_p \equiv 1 \mod M$ at all points $\pnt \ra \pntt$ 
and $f_p \equiv 1 \mod M$ at all points of 
${\Sppp{\sD} \cap \Cspt} \setminus \nu_{\Csp}^{-1}(p)$. 
Then 
\[ 
\sum_{\pnt \ra \pntt} \lrpn{\rmp{\sD}\big|_{\Cspt}, f_0}_{\pnt} 
 \;=\; \sum_{\pnt \,\in\, \Sppp{\sD} \,\cap\, \Cspt} 
        \lrpn{\rmp{\sD}\big|_{\Cspt}, f_p}_{\pnt} 
 \;=\; - \sum_{c \,\in\, \Cspt \setminus \Sppp{\sD}} \val_c(f_p) \, \rmp{\sD}(c) 
\] 
Let $T \subset \Base$ be an affine neighbourhood of $\base\pn{\Csp}$, 
and let $\Fsp \in \sO_{\CrX_T, \pi \Sppp{\sD}_{\CrX}}^*$ 
be a lift of $f_p \in \sO_{\Csp,\pntt}^*$. 
We consider $\Fsp$ as an element of $\sO_{\CrY_T, \Sppp{\sD}_{\CrY}}^*$. 
Let $\Sppp{\Fsp}$ denote the support of $\dv(\Fsp)$ in $\CrY_T$. 
Shrinking $T$ if necessary, 
we may assume that $\Sppp{\Fsp} \cap \Sppp{\sD}_{\CrY} = \varnothing$, 
that $\Sppp{\Fsp}$ intersects the fibres $\CrY_b$ over $T$ transversally 
and that $\Sppp{\Fsp} \lra T$ is \'etale. 
We choose $\Fsp$ in such a way that $\base\pn{\Csp} \in T$. 
Furthermore we may assume that all fibres $\CrY_b$ over $T$ are normal 
by Bertini's Theorem (cf.\ Proposition \ref{BertiniThm}). 
We obtain a regular map $\lam: \Sppp{\Fsp} \lra \Gp\pn{\sD}$ 
defined by 
\[ 
\lam(s) 
 \;=\; - \val_s \pn{\Fsp|_{\CrY_{\base(s)}}} \, \rmp{\sD}(s) 
 \;=\; - \val_{\Es(s)} \pn{\Fsp} \, \rmp{\sD}(s) 
\] 
for $s \in \Sppp{\Fsp}$, 
where $\Es(s)$ is the unique irreducible component of $\Sppp{\Fsp}$ 
containing $s$. 
While $\Sppp{\Fsp} \lra T$ is finite \'etale, taking the trace of $\lam$ over $T$ 
yields the map $\Trace \lam : T \lra \Lin_R$, given by 
\begin{eqnarray*}
\Trace \lam (b) 
 & = & \sum_{s \ra b} \lam(s) 
  =  - \sum_{c \,\in\, \CrY_{b} \setminus \Sppp{\sD}} 
              \val_c \pn{\Fsp|_{\CrY_{b}}} \, \rmp{\sD}(c) 
  =  - \sum_{c \,\in\, \CrY_b \setminus \Sppp{\sD}} 
              \lrpn{\rmp{\sD}|_{\CrY_b}, \Fsp|_{\CrY_b}}_{c} \\ 
 & = & \sum_{\pnt \,\in\, \CrY_b \,\cap\, \Sppp{\sD}} 
            \lrpn{\rmp{\sD}|_{\CrY_b}, \Fsp|_{\CrY_b}}_{\pnt} 
  =  \bigpn{\pn{\nu_{\CrX_b,*} \circ \fml} 
           \pn{\sD\cut\CrY_b}}_{\pi \Sppp{\sD} \,\cap\, \CrX_b} 
           \pn{\Fsp|_{\CrY_b}} 
\end{eqnarray*} 
Then $\Trace \lam \bigpn{\base\pn{\Csp}} = 
\bigpn{\pn{\nu_{\Csp,*} \circ \fml} \pn{\sD\cut\Cspt}}_{\pntt} \pn{\fsp}$. 
As $\lam$ is regular on $\Sppp{\Fsp}$, $\Trace \lam$ is regular on $T$, 
according to \cite[III, No.~5, Prop.~8]{S}. 
Thus $\Phi := \Trace \lam$ gives the desired map. 
\ePf 

\bCor 
\label{upperSemiCont}
Let $\fmlE \cong \Zint^t \tms_k \pn{\Gac}^v$ 
be a formal subgroup of $\Divf^0_{\Xt}$ 
which contains $\Divf^0_{\Xt/X}$, 
and let $\Base$ be a variety parametrizing Cartier curves in $X$. 
The function 
\bDpl \dim \lrpn{ \fmlE \tms_{\Divf_{\wt{\lul}}} \Divf_{\wt{\lul}/\lul}^0 }^{\vee} : 
   \Base \lra \Nat
\eDpl 
is upper semi-continuous. 
\eCor 

\bPf 
For $C \in \Base$ write 
\;$\dime(C) := \dim \bigpn{ \fmlE \tms_{\Divf_{\Ct}} \Divf_{\Ct/C}^0 }^{\vee}$. 
%
We have to show that the sets 
$\st{C \in \Base \:|\: \dime(C) \leq n}$ 
are open for all $n \in \Nat$. 
Let $\Csp \in \Base$ with $\dime(\Csp) = n$. 
We show: there is an open neighbourhood $U \ni \Csp$ in $\Base$ 
such that $\dime(C) \leq n$ for all $C \in U$. 

For some $c_0, c_1 \in \Nat$ with $c_0 + c_1 = c := \dim \fmlEd - n$ 
one finds $\Zint$-linearly independent elements 
$\del_1^{(0)}, \ldots, \del_{c_0}^{(0)} \in 
\pn{\fmlE \cut \Cspt}(k) \setminus \Divf_{\Cspt/\Csp}^0(k)$ 
and $k$-linearly independent elements 
$\del_1^{(1)}, \ldots, \del_{c_1}^{(1)} \in 
\Lie \pn{\fmlE \cut \Cspt} \setminus \Lie \bigpn{\Divf_{\Cspt/\Csp}^0}$ 
that extend a basis of $\Divf_{\Cspt/\Csp}^0(k)$ 
resp.\ $\Lie\bigpn{\Divf_{\Cspt/\Csp}^0}$ 
to a basis of $\fmlE(k)$ resp.\ $\Lie\pn{\fmlE}$. 
Let $\sD_j^{(i)} \in \fmlE$ with $\sD_j^{(i)} \cut \Cspt = \del_j^{(i)}$ 
for $j = 1, \ldots, c_i$ and $i = 0,1$. 
By Proposition \ref{openProperty} there exists an open neighbourhood 
$U \ni \Csp$ in $T$ such that we have 
$\bigpn{\nu_{C,*} \circ \fml} \bigpn{\sD_j^{(i)} \cut \Ct} \neq 0$ 
for $j = 1, \ldots, c_i$ and $i = 0,1$ for all $C \in U$, 
and the locus in $\Base$ where 
$\bigpn{\nu_{C,*} \circ \fml} \bigpn{\sD_1^{(i)} \cut \Ct}, \ldots, 
\bigpn{\nu_{C,*} \circ \fml} \bigpn{\sD_{c_i}^{(i)} \cut \Ct}$ 
are linearly dependent mod $\Divf_{\Ct/C}^0$ 
is a closed proper subset $V \subset U$. 
Then $\dime(C) \leq \dim \fmlEd - c = n$ for all $C \in U \setminus V$. 
\ePf 

\bPrp 
\label{FamilyShrinking}
Let $\sL$ be an ample line bundle on $X$. 
Then there exists $N\in\Nat$ such that \,$\Divf^0_{\Xt/X}$ can be computed 
from curves in the parameter space $\Base:=\left|\sL^N\right|_S^{d-1}$: 
this means by definition 
\[ \Divf_{\Xt/X}^0 = 
   \bigcap_{C\in\Base} \lrpn{\llul\cut\Ct}^{-1} \Divf_{\Ct/C}^0 \laurin 
\] 
\ePrp 

\bPf 
Cartier curves in $X$ are complete intersections locally in a neighbourhood 
of the singular locus $S$. 
For the computation of $\Divf_{\Xt/X}^0$ only the formal neighbourhood 
of $S$ is relevant, 
since $\Divf_{\Xt/X}^0$ is defined via the push-forward of formal divisors 
($\see$Definition~\ref{Div_Y/X}) 
and its support is contained in $S_{\Xt} := S \tms_X \Xt$. 
Thus we can replace the range of all Cartier curves 
by a set of complete intersection curves. 
As $\Divf_{\Xt/X}^0$ is dual-algebraic, 
there exists an effective divisor $\Es$ on $\Xt$ with support in $S_{\Xt}$ 
such that $\Divf_{\Xt/X} \subset \Fm{\Xt}{\Es}$ 
($\seecite$\cite[Prop.~3.21]{Ru2} 
or proof of \cite[Prop.~3.24]{Ru}), 
where $\Fm{\Xt}{E}$ is the formal group associated with the modulus $\Es$ 
($\seecite$\cite[Def.~3.13]{Ru2}). 
If $m$ is the maximum multiplicity of $\Es$, 
then only the $(m-1)^{\mathrm{th}}$-infinitesimal neighbourhood of $S$ 
is relevant. 
(This follows e.g.\ from \cite[Lem.~3.21]{Ru}.) 
Any Cartier curve $C \subset X$ 
can be approximated by a curve $C_N \in |\sL^N|$ 
such that the $(m-1)^{\mathrm{th}}$-infinitesimal neighbourhood of $S \cut C$ 
coincides with the one of $S \cut C_N$ 
for sufficiently large $N \in \Nat$. 
\ePf 

\bPrp 
\label{finiteCurves}
Let $\fmlE \cong \Zint^t \tms_k \pn{\Gac}^v$ 
be a formal subgroup of $\Divf^0_{\Xt}$ 
containing $\Divf^0_{\Xt/X}$, 
and let $\Base$ a parameter space of curves such that $\Divf^0_{\Xt/X}$ 
can be computed from $\Base$. 
Then there are finitely many curves $C_1,\ldots,C_r \in \Base$ 
such that 
\[ \Divf^0_{\Xt/X} = \bigcap_{i=1}^{r} 
   \left.\lrpn{\llul\cut\Cit}\right|_{\fmlE}^{-1} \Divf^0_{\Cit/C_i} \laurin 
\] 
\ePrp 

\bPf 
For $C \in \Base$ set 
$\fmlG_C := \bigpn{\llul \cut \Ct}\big|_{\fmlE}^{-1} \Divf_{\Ct/C}^0$. 
Then it holds $\Divf_{\Xt/X}^0 = \bigcap_{C \in \Base} \fmlG_C$. 
For every sequence $\st{C_{\nu}}$ of curves in $\Base$ 
the sequence $\st{\fmlE_{\nu}}$ with 
$\fmlE_0 := \fmlE$, $\fmlE_{\nu + 1} := \fmlE_{\nu} \cap \fmlG_{C_{\nu}}$ 
becomes stationary, 
due to the noetherian properties of the formal group $\fmlE$, 
cf.\ \cite[Remark 3.26]{Ru}. 
\ePf 

\bCor 
\label{oneCurve}
If $\Ci = H^{(i)}_1 \isec \ldots \isec H^{(i)}_{d-1}$, then 
\bDpl \Csp := 
   \bigpn{\sum_iH^{(i)}_1} \isec \ldots \isec \bigpn{\sum_iH^{(i)}_{d-1}}
\eDpl 
satisfies: 
\[ \Divf^0_{\Xt/X} = 
   \left.\lrpn{\llul\cut\Cspt}\right|_{\fmlE}^{-1} \Divf^0_{\Cspt/\Csp} \laurin 
\] 
\eCor 

\bThm 
\label{genCurve}
Let $\fmlE \cong \Zint^t \tms_k \pn{\Gac}^v$ 
be a formal subgroup of $\Divf^0_{\Xt}$ 
such that $\Divf^0_{\Xt/X} \subset \fmlE$. 
Let $\sL$ be an ample line bundle on $X$. 
For $C\in\left|\sL^N\right|_S^{d-1}$ and sufficiently large $N$  
the property 
\[ \Divf^0_{\Xt/X} = \left.\lrpn{\llul\cut\Ct}\right|_{\fmlE}^{-1} \Divf^0_{\Ct/C} 
\] 
is open and dense. 
\eThm 

\bPf 
By Propositions \ref{FamilyShrinking}, \ref{finiteCurves} 
and Corollary \ref{oneCurve}   there exist numbers 
$r,N\in\Nat$ and a curve 
$\Csp\in\left|\sL^{rN}\right|_S^{d-1}$ with \;
$ \Divf^0_{\Xt/X} = \big.\big(\llul\cut\Cspt\big)\big|_{\fmlE}^{-1} 
                  \Divf^0_{\Cspt/\Csp} 
$. 
Let $C$ be a general curve in $\left|\sL^{rN}\right|_S^{d-1}$. 
By definition of $\Divf^0_{\Xt/X}$ we have 
\[ \left.\lrpn{\llul\cut\Ct}\right|_{\fmlE}^{-1} \Divf^0_{\Ct/C} \;\supset\; 
   \Divf^0_{\Xt/X} \;=\; 
   \left.\lrpn{\llul\cut\Cspt}\right|_{\fmlE}^{-1} \Divf^0_{\Cspt/\Csp} \laurink 
\] 
i.e.\ 
\bDpl \dim \lrpn{ \fmlE \tms_{\Divf_{\Ct}} \Divf_{\Ct/C}^0 }^{\vee} \;\geq\; 
   \dim \lrpn{ \fmlE \tms_{\Divf_{\Cspt}} \Divf_{\Cspt/\Csp}^0 }^{\vee} \laurin 
\eDpl 
According to Pro-position \ref{upperSemiCont} 
the expression 
\;$\dim \bigpn{ \fmlE \tms_{\Divf_{\Ct}} \Divf_{\Ct/C}^0 }^{\vee}$, 
as a function in $C \in \Base$, is upper semi-continuous. 
Thus a general curve $C \in \Base$ satisfies 
\[ \left.\lrpn{\llul\cut\Ct}\right|_{\fmlE}^{-1} \Divf^0_{\Ct/C} \;=\; 
   \left.\lrpn{\llul\cut\Cspt}\right|_{\fmlE}^{-1} \Divf^0_{\Cspt/\Csp} \;=\;
   \Divf^0_{\Xt/X} \laurin 
\] 
\hfill 
\ePf


\section{Example: Product of two Cuspidal Curves} 
\label{sec:Exm_Surf}

We conclude this paper with the discussion of an example that 
was the subject of the diploma of Alexander Schwarzhaupt \cite{Sch}. 
This example illustrates some pathological properties: 
the Albanese of Esnault-Srinivas-Viehweg 
is not in general compatible with products, 
in this example we obtain 
(writing $\Albsg{X} := \Urq{X}$) 
\[ \dim\big(\Albsg{\Galp\tms\Gbet}\big) > 
   \dim\big(\Albsg{\Galp}\tms\Albsg{\Gbet}\big) \laurin 
\]  
Moreover, given a very ample line bundle $\sL$ 
on the surface $X=\Galp\tms\Gbet$ 
and a curve $C_N \in \left|\sL^N\right|$ in general position, 
we work out a necessary and sufficient condition on the integer $N$ 
for the surjectivity of the Gysin map $\Albsg{ C_N } \lra \Albsg{ X }$. 

\noindent 
The base field $k$ is assumed to be algebraically closed and of characteristic 0.

\subsection{Cuspidal Curve} 
\label{sec:Galp}

Let $\Galp \subset \Prj^2_k$ be the projective curve defined by 
\[ \Galp: \quad X^{2\alp +1} - Y^2 Z^{2\alp -1}=0 
\] 
where $X:Y:Z$ are homogeneous coordinates of $\Prj^2_k$ and 
$\alp \geq 1$ is an integer. The singularities of this curve are cusps 
at $0:=[0:0:1]$ and $\infty:=[0:1:0]$. 
The normalization $\Galpt$ of $\Galp$ is the projective line: 
\bDpl \Galpt = \Prj^1_k \laurin 
\eDpl 
Then $\Alb\pn{ \Galpt } = \Alb\pn{ \Prj^1_k } = 0$. Since 
$\Albsg{\Galp}$  is an extension of $\Alb\pn{ \Galpt }$ 
by the linear group $L_{\Galp} = \bigpn{ \Divf^0_{\Galpt / \Galp}}^{\vee}$, 
we obtain 
\[ \Albsg{ \Galp } = L_{\Galp} = \lrpn{ \Divf^0_{\Galpt / \Galp}}^{\vee} \laurin 
\] 
Moreover, $\Galp$ is homeomorphic to $\Prj^1_k$, i.e.\ the normalization 
$\Galpt$ is given by the largest homeomorphic curve $\Galp'$. 
This implies that $L_{\Galp}$ is a unipotent group 
($\see$Theorem~\ref{cusp}) 
and equivalently 
$\Divf^0_{\Galpt / \Galp}$ is an infinitesimal formal group 
($\seecite$\cite[Proposition 1.17]{Ru2}). 
We have 
\bDpl \Lie \lrpn{ \Divf^0_{\Galpt/\Galp} } = \Homk \bigpn{ L_{\Galp}(k),k } \laurin 
\eDpl 

The $k$-valued points of $L_{\Galp}$ are given by ($\see$Theorem~\ref{cusp}) 
\[  L_{\Galp}(k) 
    \;=\;  \frac{1+\fm_{\Galpt,0}}{1+\fm_{\Galp,0}} \tms 
               \frac{1+\fm_{\Galpt,\infty}}{1+\fm_{\Galp,\infty}}  
    \;\cong\;  \frac{\fm_{\Galpt,0}}{\fm_{\Galp,0}} \oplus 
                      \frac{\fm_{\Galpt,\infty}}{\fm_{\Galp,\infty}} \laurin 
\] 
The dimensions are computed in \cite[Proposition 1.5]{Sch} as 
\begin{eqnarray*} 
\dim_k \bigpn{ \fm_{\Galpt,0} \big/ \fm_{\Galp,0} } & = & \alp \\
\dim_k \bigpn{ \fm_{\Galpt,\infty} \big/ \fm_{\Galp,\infty} } & = & 
          2 \alp (\alp -1)
\end{eqnarray*} 
hence 

\bPrp 
\label{dimAlbGam}
\[ \dim\,\Albsg{ \Galp }  =  \dim\, L_{\Galp}  
      =  \dim_k \,\Lie \lrpn{ \Divf^0_{\Galpt/\Galp} }  
      =  \alp (2 \alp -1) \laurin 
\] 
\ePrp 

As 
\;$\Lie \bigpn{\Divf^0_{\Galpt/\Galp}} 
 = \fml^{-1} \dsum_{q=0,\infty} \Hom_k \bigpn{\fm_{\Galpt,q} \big/ \fm_{\Galp,q}, k}$\; 
it holds 
\[ \# \, \val_q \lrpn{ \Lie \bigpn{ \Divf^0_{\Galpt/\Galp} } } 
 =  \# \lrpn{ \val_q (\sO_{\Galpt}) \setminus \val_q (\sO_{\Galp}) } 
 =  \dim_k \bigpn{ \fm_{\Galpt,q} \big/ \fm_{\Galp,q} } \laurink 
\] 
for $q \in \{0,\infty \} \subset \Galpt$, cf.\ \cite[No.~3.3]{Ru}. 
A basis of \;$\Lie \bigpn{ \Divf^0_{\Galpt/\Galp} }$ is given by 
the following set of representatives 
\[  \Theta_{\Galp} = \left\{ \, t_q^{-\nu} \;\left| \; -\nu \in 
    \val_q\lrpn{\Lie\bigpn{\Divf^0_{\Galpt/\Galp}}},\; q=0,\infty\right.\right\}
\]
where $t_q$ is a local parameter of $\fm_{\Galpt,q}$ at $q$ 
and $t_q^{-\nu} \in \sO_{\Galpt,p}$ for all $p \neq q$.

\subsection{Cuspidal Surface} 
\label{sec:Galp_x_Gbet}

Let $X$ be the product of the cuspidal curves $\Galp,\Gbet$ 
from Subsection~\ref{sec:Galp}: 
\[ X = \Galp \times \Gbet \] 
where $\alp,\beta \geq 1$ are integers. 
The singular locus of $X$ is 
\[ X_{\sing} = \bigpn{\st{0}\tms\Gbet} \cup \bigpn{\st{\infty}\tms\Gbet} \cup 
             \bigpn{\Galp\tms \st{0}} \cup \bigpn{\Galp\tms\st{\infty}} \laurin 
\] 
The normalization $\Xt$ of $X$ is a resolution of singularities 
and given by 
\bDpl \Xt = \Galpt \times \Gbett = \Prj^1_k \times \Prj^1_k \laurin 
\eDpl 
Then 
$\Alb\pn{ \Xt } = \Alb\pn{ \Prj^1_k } \times \Alb\pn{ \Prj^1_k } = 0$. 
Thus the Albanese of Esnault-Srinivas-Viehweg of $X$ 
coincides with its affine part: 
\[ \Albsg{X} = L_X = \lrpn{ \Divf^0_{\Xt/X} } ^{\vee} \] 
and $\Divf^0_{\Xt/X}$ is an infinitesimal formal group, 
since the normalization is a homeomorphism. 
The task is now to determine $\Divf_{\Xt/X}$. 
The support of $\Divf^0_{\Xt/X}$ is the preimage of $X_{\sing}$: 
\[ \Supp \lrpn{ \Divf^0_{\Xt/X} } = \bigpn{\st{0}\tms\Gbett} 
   \cup \bigpn{\st{\infty}\tms\Gbett} 
   \cup \bigpn{\Galpt\tms \st{0}} \cup \bigpn{\Galpt\tms\st{\infty}} \laurin 
\] 
By Proposition~\ref{genCurve} we can compute $\Divf_{\Xt/X}^0$ 
by the preimage of $\Divf_{\Ct/C}^0$ in $\Divf_{\Xt}$ 
under pull-back 
to $\Ct$, for a sufficiently ample curve $C \subset X$ in general position. 
We may take 
$C = \bigcup_{i=1}^k \bigpn{\st{p_i} \tms \Gbet} 
\cup \bigcup_{j=1}^l \bigpn{\Galp \tms \st{q_j}}$ 
for sufficiently many points $p_i \in \Galp\setminus\st{0,\infty}$ 
and $q_j \in \Gbet\setminus\st{0,\infty}$. 
Then $\Divf_{\Ct/C}^0 = 
\prod_{i=1}^k \Divf_{\Gbett/\Gbet}^0 \tms \prod_{j=1}^l \Divf_{\Galpt/\Galp}^0$. 
Now \;$\Lie \bigpn{ \Divf^0_{\Xt/X} }$\; is the space of those 
$\del \in \Lie \bigpn{ \Divf^0_{\Xt} }$ with support in the preimage of $X_{\sing}$ 
such that $\del \cut \bigpn{\st{p} \tms \Gbett} \in \Lie\bigpn{\Divf_{\Gbett/\Gbet}^0}$ 
and $\del \cut \bigpn{\Galpt \tms \st{q}} \in \Lie\bigpn{\Divf_{\Galpt/\Galp}^0}$ 
for all $p \in \Galpt\setminus\st{0,\infty}$ and all $q \in \Gbett\setminus\st{0,\infty}$. 
From this we see: 
If $\Theta_{\Gamma_{\iota}}$ is a set of representatives of a basis of 
$\Lie \bigpn{ \Divf^0_{\widetilde{\Gamma}_{\iota} / \Gamma_{\iota}} }$ 
for $\iota=\alp,\beta$, then
\[ \Theta_{\Galp\tms\Gbet} = 
   \Bigst{ \, \tha_{\Galp} \tens \tha_{\Gbet} \;\Big| \; 
   (\tha_{\Galp},\tha_{\Gbet}) \in 
   \Bigpn{ \bigpn{\Theta_{\Galp} \cup \{1\}} \tms 
                 \bigpn{\Theta_{\Gbet} \cup \{1\}} } 
   \setminus \{(1,1)\} }
\]
is a set of representatives of a basis of \;$\Lie \bigpn{ \Divf^0_{\Xt/X} }$. 

Thus the dimension of $\Albsg{X}$ is given by 
\begin{eqnarray*}
&   & \dim\, \Albsg{ \Galp \tms \Gbet }  
\; = \; \dim_k \, \Lie\lrpn{\Divf^0_{\pn{\Galpt \tms \Gbett}/\pn{\Galp \tms \Gbet}}}  \\ 
& = & \lrpn{ \dim_k \, \Lie\bigpn{\Divf^0_{\Galpt/\Galp}} +1 } \cdot 
          \lrpn{ \dim_k \, \Lie\bigpn{\Divf^0_{\Gbett/\Gbet}} +1 } -1     \\
& = & \bigpn{\alp (2 \alp -1) +1} \cdot \bigpn{\beta (2 \beta -1) +1} -1 \laurin 
\end{eqnarray*} 
With Proposition~\ref{dimAlbGam} this yields 

\bPrp 
\label{dimAlb(CxC)}
\[ \dim \Albsg{ \Galp \tms \Gbet }  
= \bigpn{ \dim \Albsg{ \Galp } +1 } \cdot \bigpn{ \dim \Albsg{ \Gbet } +1 } - 1 \laurin 
\] 
\ePrp 

We obtain a basis of \;$ \im \lrpn{ \Lie\bigpn{\Divf^0_{\Xt/X}} \lra 
  \lrpn{ \left. \sK_{\Xt} \right/ \sO_{\Xt} }_{\Galpt \tms \st{q}} } $    \\
for $q \in \{0,\infty \} \subset \Gbett$ from the following set of representatives 
\[ \Theta_{\Galp\tms q} = \left\{ \tha_{\Galp} \tens t_q^{-\nu} \;\left| 
   \; \tha_{\Galp} \in (\Theta_{\Galp} \cup \{1\}), \; 
   -\nu \in \val_q \lrpn{ \Lie\bigpn{\Divf^0_{\Gbett/\Gbet}} } \right. \right\} \laurin 
\] 
Now \;$\val_q \lrpn{ \Lie\bigpn{\Divf^0_{\Gbett/\Gbet}} } = 
     \val_{\Galpt \tms \st{q}} \lrpn{ \Lie\bigpn{\Divf^0_{\Xt/X}} } $, therefore 

\bPrp 
\label{dimLie_Cxq}
\[ \dim_k \, \lrpn{ \Lie\bigpn{\Divf^0_{\Xt/X}} }_{\Galpt \tms \st{q}}^{-\nu} 
   \;=\; \dim_k \, \Lie\bigpn{\Divf^0_{\Galpt/\Galp}} + 1 
\] 
for $q \in \{0,\infty \} \subset \Gbett$ 
and all $-\nu \in \val_{\Galpt \tms \st{q}} \bigpn{ \Lie\bigpn{\Divf^0_{\Xt/X}} }$, 
where we use the notation from Lemma~\ref{.C inj}. 
Analogously for $\st{p} \times \Gbett$, $p \in \{0,\infty \} \subset \Galpt$. 
\ePrp

\subsection{Gysin Map} 
\label{sec:map:C_N->X}

Consider the following divisor $D^{k,l}$ on $X$ 
\[ D^{k,l} = \sum_{i=1}^k \bigpn{\st{p_i} \tms \Gbet} + 
                   \sum_{j=1}^l \bigpn{\Galp \tms \st{q_j}} 
\] 
where $p_i \in \Galp\setminus\st{0,\infty}$ for $i=1,\ldots,k$ 
and $q_j \in \Gbet\setminus\st{0,\infty}$ for $j=1,\ldots,l$. 
The normalization of $D^{k,l}$ is isomorphic to the disjoint union of $k+l$ 
copies of $\Prj^1_k$. 
Therefore the Picard variety of the normalization $\Pic^0\wt{D^{k,l}}$ is trivial. 
Then by Theorem~\ref{Main}, using the explicit formulas of  Propositions~\ref{multi pt} and \ref{cusp}, 
\[ \Pic^0 D^{k,l} = \Trs \times \Vcl 
\]
where $\Trs \cong (\Gm)^t$ is a torus of rank 
\[ t  =  \# S_2 - \# \Cp \lrpn{ D^{k,l} } + 1 
   =  k \cdot l - (k+l) + 1 
   =  (k-1) \cdot (l-1) \laurink 
\] 
and $\Vcl \cong (\Ga)^v$ is a vectorial group of dimension 
\[ v  =  k \cdot \dim \, L_{\Gbet} + l \cdot \dim \, L_{\Galp}  
   =  k \cdot \beta (2 \beta -1) + l \cdot \alp (2 \alp -1) \laurin 
\] 
For general $p_i \in \Galp$ and $q_j \in \Gbet$ the divisor $D^{2\alp +1,2\beta +1}$ 
is very ample ($\see$\cite{Sch} Lemma 3.2). Set 
$\sL = \sO \lrpn{ D^{2\alp +1,2\beta +1} } $ and choose 
$C_N \in |\sL^N|$ in general position. 
As $\dim \Pic^0 C_N = \textrm{const.}$ 
and $\dim \Pic^0 \bigpn{C_N \tms_X \Xt} = \textrm{const.}$ 
among $C_N \in \big|\sL^N\big|$ 
and $C_N \tms_X \Xt = {C_N}'$ is the semi-normalization, 
the dimension of 
the vectorial part $\Vcl_{C_N} = \ker \bigpn{\Pic^0 C_N \lra \Pic^0 {C_N}' }$ 
of $\Pic^0 C_N = \Albsg{ C_N } $ is constant among $C_N \in |\sL^N|$ 
and hence 
\[ \dim \Vcl_{C_N} 
   = \dim \Vcl_{N \, D^{2\alp +1,2\beta +1}} 
   = N \big( \alp (2 \alp -1)(2\beta+1) + \beta (2 \beta -1)(2\alp+1) \big) \laurin 
\] 
Since $\Albsg{\Galp\tms\Gbet}$ is a vectorial group, 
the map \\ 
\;$\Albsg{ C_N } \lra \Albsg{ \Galp \tms \Gbet }$\; 
cannot be surjective if \\ $\dim \Vcl_{C_N} < \dim\, \Albsg{ \Galp \tms \Gbet }$. \\ 
 Therefore a comparison of dimensions yields: 

\bPrp 
\label{low_bound}
The Gysin map $\Albsg{ C_N } \lra \Albsg{ \Galp \tms \Gbet }$ is not  
surjective for 
\[ N < \frac{(\alp (2 \alp -1) +1) \cdot (\beta (2 \beta -1) +1) -1} 
            { \alp (2 \alp -1)(2\beta+1) + \beta (2 \beta -1)(2\alp+1)}
\]
In the case $\alp = \beta$, this expression simplifies to 
\[ N < \frac{\alp (2\alp-1)+2}{2(2\alp+1)} \laurin 
\]
\ePrp 

The homomorphism of vectorial groups $\Vcl_{C_N} \lra \Vcl_X = \Albsg{X}$ 
is dual to the map between Lie algebras 
$\lul\cut\,\Ct_N:\Lie \bigpn{ \Divf^0_{\Xt/X} } \lra \Lie \bigpn{ \Divf^0_{\Ct_N /C_N} }$, 
and the surjectivity of the first homomorphism is equivalent to the 
injectivity of the latter one. 
Here Definition \ref{Div_Y/X} of $\Divf^0_{\Xt/X}$ 
implies immediately that the image of $\Divf^0_{\Xt/X}$ under 
pull-back $\lul\cut\,\Ct: \Divf^0_{\Xt/X} \lra \Divf^0_{\Ct}$ 
is contained in $\Divf^0_{\Ct/C}$. 

The estimation of Proposition~\ref{low_bound} 
yields a necessary condition for surjectivity of the Gysin map, 
i.e.\ a bound for $N$ from below. 

The criterion of Lemma~\ref{.C inj} 
gives a sufficient condition for the surjectivity of the Gysin map: 
\bDpl \# \lrpn{ \Ct_N \cap \bigpn{\Galpt \tms \st{q}} } 
   \; \geq \; \dim_k \, \lrpn{ \Lie\bigpn{\Divf^0_{\Xt/X}} }_{\Galpt \tms \st{q}}^{-\nu} 
\eDpl  \\ 
for all $-\nu \in \val_{\Galpt \tms \st{q}} (\Lie \Divf^0_{\Xt / X})$, 
all $q \in \{0,\infty \} \subset \Gbett$, 
and the same formula with $\bigpn{\Galpt \tms \st{q}}$ replaced by 
$\bigpn{\st{p} \tms \Gbett}$ for all $p \in \{0,\infty \} \subset \Galpt$. 
Since 
\;$\dim_k  \bigpn{ \Lie\bigpn{\Divf^0_{\Xt/X}} }_{\Galpt \tms \st{q}}^{-\nu} 
 \,=\, \dim_k \, \Lie\bigpn{\Divf^0_{\Galpt/\Galp}} + 1 $\; 
($\see$Proposition~\ref{dimLie_Cxq}) 
this is equivalent to 
\bDpl \# \lrpn{ \Ct_N \cap \bigpn{\Galpt \tms \st{q}} } \,\geq\, 
   \dim_k \, \Lie \lrpn{ \Divf^0_{\Galpt/\Galp} } + 1 \laurin 
\eDpl \\ 
Then since 
\bDpl   \# \lrpn{ \Ct_N \cap \bigpn{\Galpt \tms \st{q}} } 
    =  N \deg \lrpn{D^{2\alp +1,2\beta +1}}_{\Xt}  
    =  N (2 \alp + 1 + 2 \beta + 1) 
    =  N \,2(\alp + \beta + 1)  
    \laurink 
\eDpl 
where 
\bDpl \lrpn{D^{k,l}}_{\Xt} = 
   \sum_{i=1}^k \bigpn{\st{p_i} \tms \Gbett} + \sum_{j=1}^l \bigpn{\Galpt \tms \st{q_j}} 
   \laurink 
\eDpl \\ 
we obtain with Proposition \ref{dimAlbGam} 

\bPrp 
\label{Gysin surjective}
The Gysin map $\Albsg{ C_N } \lra \Albsg{ \Galp \tms \Gbet }$ is surjective if 
\[ N \geq \frac{\alp(2\alp-1)+1}{2(\alp+\beta+1)} 
   \qquad \textrm{and} \qquad 
   N \geq \frac{\beta(2\beta-1)+1}{2(\alp+\beta+1)} \laurin 
\]
For $\alp = \bet$ this condition is necessary and sufficient. 
\ePrp 

\bPf 
The first statement follows from the discussion above. 
In the case $\alp = \bet$ we need to show 
that the estimation above 
and the formula from Proposition~\ref{low_bound} 
yield the same bound for $N \in \Nat$. 
As the difference is $\frac{1}{2(2 \alp + 1)} < 1$, 
it suffices to check that 
$2 (2 \alp + 1)$ does not divide 
$\alp(2\alp - 1) + 1 = \alp(2\alp + 1) - (2\alp - 1)$, 
which is obvious. 
\ePf 

\vspace{\vs}

In \cite[Variant 6.4]{ESV} the following sufficient condition for surjectivity 
of the Gysin map is given: 
\[ \dim_k \, \im \Bigpn{ \H^0 \lrpn{ X,\sL^N } \lra \H^0 \lrpn{ Z,\sL^N |_Z } } 
   \,\geq\, 2 \dim L_{C_N} + \# \Cp \lrpn{ X } + 2 
\]
for all $Z \in \Cp(X)$, 
where $L_{C_N}$ is the largest connected  affine subgroup of 
$\Pic^0 C_N$ for $C_N \in \left| \sL^N \right|$ in general position. 
For $X = \Galp \tms \Gbet$ 
it holds $\Cp\pn{X} = \st{X}$ and $\Pic^0 C_N = L_{C_N} = \Vcl_{C_N}$. 
Alexander Schwarzhaupt showed in his diploma \cite{Sch} 
that in our example and for $\alp = \beta$ 
this condition leads to the estimation 
\[ N > 2 \frac{3 \alp (2 \alp - 1) - 1}{2 \alp + 1} + 1 \laurin 
\]

\end{document}